\documentclass[letterpaper, 10 pt, conference]{ieeeconf}  % Comment this line out
\IEEEoverridecommandlockouts                              % This command is only
\usepackage{braket,amsfonts,multirow,graphics,siunitx}
\usepackage{latexsym,amsmath,amssymb,textcomp}
\allowdisplaybreaks
\usepackage{balance}
\usepackage{graphicx}
\usepackage{epstopdf}
\usepackage{tabularx}
\usepackage{booktabs}
\usepackage{subcaption}
\usepackage{lipsum}
\usepackage{pgfplots}
\usepackage{xcolor}
\usepackage{caption}
\usepackage{pstool}
\usepackage{tikz}
\usepackage{amsmath}

\usetikzlibrary{shapes.geometric, arrows}

\tikzstyle{block} = [draw, rectangle, minimum height=.7cm, minimum width=.7cm]
\tikzstyle{sum} = [draw, circle, inner sep=1pt]
\tikzstyle{input} = [coordinate]
\tikzstyle{output} = [coordinate]
\tikzstyle{arrow} = [thick, ->, >=stealth]
\newtheorem{thrm}{Theorem}[section]

\newtheorem{lemma}[thrm]{Lemma}

\newtheorem{cor}[thrm]{Corollary}
\newtheorem{proposition}[thrm]{Proposition}

\newtheorem{rmrk}[thrm]{Remark}

\newcommand{\R}{{\mathbb R}}

\newcommand{\Real}{\mathbb{R}}

\title{System Realizations by Mammillary Models\\%
with an Application to Propofol Pharmacokinetics}%
\author{Veronica Beltrami$^1$, Luca Consolini$^1$, Mattia Laurini$^1$, Marco Milanesi$^{2}$, Michele Schiavo$^{2}$, Antonio Visioli$^{2}$%
\thanks{This work was supported by the PRIN 2022 Project ``ACTIVA -- Automatic Control of Total IntraVenous Anesthesia'',
funded by European Union -- Next Generation EU under Grant CUP D53D23001180006.}
\thanks{$^1$ Department of Engineering and Architecture,
University of Parma, 43124 Parma, Italy,
e-mail: \{{\tt\scriptsize veronica.beltrami}, {\tt\scriptsize luca.consolini}, {\tt\scriptsize mattia.laurini}\}{\tt\scriptsize @unipr.it}}
\thanks{$^{2}$ Department of Mechanical and Industrial Engineering,
University of Brescia, 25136 Brescia, Italy,
e-mail: \{{\tt\scriptsize marco.milanesi}, {\tt\scriptsize michele.schiavo}, {\tt\scriptsize antonio.visioli}\}{\tt\scriptsize @unibs.it}}%
}
\date{}

\begin{document}

\maketitle

\begin{abstract}
This work addresses the problem of linear system realizations by mammillary models, offering necessary and sufficient conditions under which a given transfer function can be represented in this form. While standard identification techniques may yield transfer functions without an explicit connection to underlying physiological processes, compartmental models, particularly mammillary ones, reflect the physiological dynamics of distribution and elimination. This feature is especially relevant in clinical pharmacology, where model parameters must correspond to meaningful biological processes to support interpretability, personalization, and safe drug delivery, such as in total intravenous anesthesia. To conclude, an application to a propofol infusion model illustrates how mammillary realizations can support physiologically interpretable system representations.
\end{abstract}

%{\small {\bf Keywords}: Depth of hypnosis control, feedforward control, induction, optimization.}

\section{Introduction}
Compartmental systems are fundamental tools for modeling the dynamics of substance distribution and elimination within complex biological and pharmacological systems.
They consist of a finite number of compartments that exchange material with each others.
The transfer of material between compartments can be described by first-order ordinary differential equations, based on the principle of mass conservation. 
These models are commonly used to describe the distribution and elimination of a drug within a human body. In this context, a human body is conceptualized as a set of interconnected compartments, each representing regions or tissues of the body that behave similarly in terms of drug absorption, distribution, and elimination, with transfer rates characterizing the movement of the drug between compartments.

%These models represent a system as a set of interconnected compartments, each containing a homogeneous quantity of material.
%In the classical formulation, the substance exchange is governed by first-order ordinary differential equations~\cite{}.
Such models are often designed under the assumption of mass conservation, meaning that any increase or decrease in the amount of substance within a compartment is fully accounted for by the flows between compartments and by explicit inputs or outputs (e.g., drug administration or elimination)~\cite{anderson1983lecture}.
%Compartmental models can be represented effectively using state-space models, which offer a clear description, particularly when dealing with multi-input multi-output systems.
%Alternatively, autoregressive with exogenous input (ARX) models provide a discrete-time statistical framework for identifying system dynamics from input-output data, offering advantages in estimation and computational efficiency~\cite{}.
Within this broad framework, mammillary compartmental models have emerged as especially significant~\cite{jacquez1972compartmental}.
These models are structured around a central compartment that interacts directly with peripheral compartments, which themselves do not exchange material directly.
This configuration is highly relevant in pharmacokinetics, where mammillary models are widely employed to characterize drug absorption, distribution, and elimination, especially in clinical scenarios such as total intravenous anesthesia control~\cite{minto1997influence, schnider1998influence}.

%In this paper we consider  
%The main focus of this problem is finding if a given transfer function can be the one of a mammillary system, and, in case it is, finding such system.  biological interpretation (see, e.g.,~\cite{}).
%Although realization problems have been extensively studied in general linear systems theory, including positive realizations and minimal realizations (see, for instance,~\cite{}), the conditions ensuring compartmental or, more specifically, mammillary realizations are more restrictive and require dedicated analysis.

%In this paper, we address this problem explicitly by investigating the realization problem for a specific subclass of linear, time-invariant mammillary models characterized by an external input exclusively into the central compartment, {\color{red} with the observed output corresponding solely to the state of this central compartment}.
%We derive necessary and sufficient conditions under which an arbitrary rational transfer function admits a mammillary realization in this form.
%Moreover, we show the practical usability of our theoretical results by applying this constructive methodology to the pharmacokinetic modeling of propofol, a widely administered intravenous anesthetic in general anesthesia.

\subsection{Main Results}
Let 
\begin{equation}\label{sys:cont}
\begin{cases}
\dot{x}(t) = A x(t) + B u(t) \\
y(t) = C x(t)
\end{cases}    
\end{equation}
be a continuous-time linear system with $A = (a_{ij})\in \R^{n \times n}$, $B = (b_i) \in \mathbb{R}^{n \times 1}$, and $C = (c_j) \in \mathbb{R}^{1 \times n}$.
According to~\cite{ anderson1983lecture, benvenuti2002positive, kajiya1985compartmental, jacquez1972compartmental},~\eqref{sys:cont} is said to be a \textit{compartmental system} if, for $i, j \in \{1, \ldots, n\}$, 
\begin{align}
b_i &\ge 0,  \quad c_j \ge 0; \label{eq:cond1} \\
a_{ij} &\ge 0 \quad \text{for } i \neq j; \label{eq:cond2} \\
a_{ii} &+ \sum_{j \neq i} a_{ji} \leq 0.\label{eq:cond3}
\end{align}

The usual physical interpretation is that, for $i \in \{1, \ldots, n\}$, state variable $x_i(t)$ (the $i$th component of vector $x(t)$) represents the amount of resource present in the $i$-th compartment at time $t$.
Each diagonal element $a_{ii}$ of matrix $A$ represents the cumulative outgoing flows from the $i$th compartment, which accounts for material loss.
Conversely, each off-diagonal element $a_{ij}$, for $i \neq j$, represents a constant incoming flow into the $i$th compartment coming from the $j$th one.
Input $u(t)$ represents an external injection of material into the system and $B$ determines how such material is distributed among the compartments.
Output $y(t)$ represents the measured quantity of material in the system at time $t$, as determined by observation vector $C$.
This reflects the amount of material present in one or more compartments, depending on which components of $C$ are nonzero.

The following is a well-established property of compartmental systems (see, for instance,~\cite{anderson1983system}).
\begin{proposition}
\label{prop:eigenvalue}
Let $A$ be the matrix associated with a compartmental system (i.e., let $A$ satisfy~\eqref{eq:cond1}--\eqref{eq:cond3}), and let $\lambda$ be an eigenvalue of $A$. Then $\operatorname{Re}(\lambda) \leq 0$.
\end{proposition}

A \textit{mammillary model} is a compartmental system that consists of a central compartment connected to a set of peripheral ones.
In this model, there is no material transfer between peripheral compartments, and all peripheral compartments exchange material with the central one.
In this work, we consider mammillary models belonging to the following class:
\begin{equation}
\begin{gathered}
\label{eq:mammillary}
A = \left[
\begin{array}{c|cccc}
\scriptstyle \!\!-k_{10}-\sum\limits_{i=2}^n k_{1i}\!\!  & \scriptstyle k_{21} & \scriptstyle k_{31} & \cdots & \scriptstyle k_{n1} \\
\hline
\scriptstyle k_{12} & \scriptstyle \!-k_{21}\! & 0 & \cdots & 0 \\
\scriptstyle k_{13} & 0 & \scriptstyle \!-k_{31}\! & \cdots & 0 \\
\vdots & \vdots & \vdots & \ddots & \vdots \\
\scriptstyle k_{1n} & 0 & 0 & \cdots & \scriptstyle \!-k_{n1}\!
\end{array}
\right] \!\!\in \R^{n \times n},\vspace{3pt}\\
B = [1, 0, \ldots, 0]^T, \qquad C = B^T = [1, 0, \ldots, 0].
\end{gathered}
\end{equation}

We assume that all $k_{i1}$, with $i \in \{2,\ldots,n\}$, are distinct.

Moreover, throughout this work, we consider that the parameters $k_{i1}$, for $i \in \{ 2, \ldots, n\}$, are ordered so that $k_{21} < k_{31} < \cdots < k_{n1}.$
Note that this ordering of the parameters is assumed without loss of generality, as there always exists a permutation of the states that realizes such an order.

Observe that, if $k_{10}>0$ and $k_{ij}>0$, for all $i \neq j \in \{1, \ldots, n\}$, then, $A, B, C$ in~\eqref{eq:mammillary} satisfy conditions~\eqref{eq:cond1}--\eqref{eq:cond3}. 

The following property is known in the literature; however, for the sake of completeness, we provide a brief proof tailored to our specific case.
\begin{proposition}
\label{thrm:eigenvalue}
Let $A$ be a matrix defined as in~\eqref{eq:mammillary}, with $k_{10}>0$ and $k_{ij}>0$, for all $i \neq j \in \{1, \ldots, n\}$.
If $\lambda$ is an eigenvalue of $A$, then $\lambda \in \mathbb{R}$ and $\lambda \le 0$.
\end{proposition}

\begin{proof}
Since $k_{10}>0$ and $k_{ij}>0$, for all $i,j \in \{1, \ldots, n\}$, then $A$ satisfies conditions~\eqref{eq:cond2} and~\eqref{eq:cond3}, and, therefore, Proposition~\ref{prop:eigenvalue} holds.
So we need only to prove that the eigenvalues of $A$ are real.

Let $D$ be a diagonal matrix of the following form
\[
D = \begin{bmatrix}
1      & 0      & 0      & \cdots & 0 \\
0      & \sqrt{k_{21}/k_{12}} & 0      & \cdots & 0 \\
0      & 0      & \sqrt{k_{31}/k_{13}} & \cdots & 0 \\
\vdots & \vdots & \vdots & \ddots & \vdots \\
0      & 0      & 0      & \cdots & \sqrt{k_{n1}/k_{1n}}
\end{bmatrix}.
\]

It is easy to compute $D A D^{-1} =: \tilde{A}$
\begin{equation}
\label{eq:A_tilde}
\tilde{A} = \left[
\begin{array}{c|cccc}
\scriptstyle \!\!-k_{10} - \sum\limits_{i=2}^n k_{1i}\!\! & \scriptstyle \!\sqrt{k_{21}k_{12}}\! & \scriptstyle \!\sqrt{k_{31}k_{13}}\!\! & \cdots & \!\!\scriptstyle \sqrt{k_{n1}k_{1n}} \\
\hline
\scriptstyle \sqrt{k_{21}k_{12}} & \scriptstyle -k_{21} & 0 & \cdots & 0 \\
\scriptstyle \sqrt{k_{31}k_{13}} & 0 & \scriptstyle -k_{31} & \cdots & 0 \\
\vdots & \vdots & \vdots & \ddots & \vdots \\
\scriptstyle \sqrt{k_{n1}k_{1n}} & 0 & 0 & \cdots & \scriptstyle -k_{n1}
\end{array}
\right].
\end{equation}
Matrices $\tilde{A}$ and $A$ are similar, therefore they have the same eigenvalues.
Moreover, $\tilde{A}$ is symmetric, so all eigenvalues of $\tilde{A}$, and hence of $A$, are real.
\end{proof}

In the first part of this paper, we investigate a realization problem.
Our objective is to derive necessary and sufficient conditions under which a given rational transfer function
$$H(s):=\frac{\beta(s)}{\alpha(s)} \text{ of order } n,$$
 with $\alpha(s)=s^n+\alpha_{n-1}s^{n-1}+\alpha_{n-2}s^{n-2}+\ldots+\alpha_1 s+ \alpha_0$,  $\alpha_i \in \R$, admits a realization of the form
\[
H(s) = C(sI - A)^{-1}B,
\]
where $I$ denotes the $n \times n$ identity matrix, and $A$, $B$ and $C$ have the specific structure described in~\eqref{eq:mammillary}. 

Recall that, in general,
\[
C (s I -A)^{-1} B=\frac{C(s I -A)^{\text{adj}}B}{\chi(s)},
\]
where $\chi(s)$ denotes the characteristic polynomial of $A$, that is $\chi(s)=\text{det}(sI-A)$.

\begin{rmrk}
\label{remark}
Note that, if in~\eqref{eq:mammillary} it holds that $k_{10}>0$ and $k_{ij}>0$, for all $i,j \in \{1, \ldots, n\}$, then it follows from Proposition~\ref{thrm:eigenvalue} that all roots of $\chi(s)$ are real and non-negative.
\end{rmrk}

The essential results are enclosed in the two following theorems.
Let $p$ be a vector collecting the parameters of $A$ in~\eqref{eq:mammillary} (i.e., the entries of $A$):
\[
p := (k_{10}, k_{21}, \ldots, k_{n1}, k_{12}, \ldots, k_{1n}) \in \mathbb{R}^{2n - 1}.
\]

\begin{thrm}
\label{thrm:conditions}
There exists a unique $p$ such that $$C(sI - A)^{-1}B=H(s)$$ if and only if the following conditions hold:
\begin{enumerate}
\item the relative degree of $H$ is 1;

\item the numerator of $H$ is a monic polynomial with simple, real, and nonzero roots.
\end{enumerate}

\end{thrm}

A consequence of Theorem \ref{thrm:conditions} is the following.
\begin{thrm}\label{thrm:conditions_p}
There exists a unique, strictly positive $p$ such that $H(s) = C(sI - A)^{-1}B$ if and only if the following conditions hold:
\begin{enumerate}
\item the relative degree of $H$ is 1;

\item the numerator of $H$ is a monic polynomial with simple, real, and strictly negative roots;

\item $H(0)>0$;

\item for any root $z_i$ of the numerator of $H$, $\frac{\beta_i(z_i)}{\alpha(z_i)}>0$.
Where $\beta_i(s):=z_i \widehat{\beta(s)}^i$ is a polynomial of degree $n-1$, and $\widehat{\beta(s)}^i$ is $\beta(s)$ without the factor $(s-z_i)$.
\end{enumerate}

\end{thrm}
As we will show in Section \ref{proof}, as a consequence of Theorem \ref{thrm:conditions}, it is possible to derive an algorithm to determine the $2n-1$ parameters

\medskip

\fbox{%
\parbox{0.95\linewidth}{%
\textit{Algorithm for computing the parameters}

Let $z_2 > z_3 > \dots > z_n$ be the roots of $\beta(s)$.

\medskip

\[
\begin{aligned}
k_{10} &= \frac{\alpha(0)}{\beta(0)} \\
k_{i1} &= -z_i, & \quad & i \in \{2,\dots,n\} \\
k_{1i} &= \frac{\alpha(z_i)}{\beta_i(z_i)}, & \quad & i \in \{2,\dots,n-1\}.
\end{aligned}
\]

\medskip

Where
\[
\beta_i(s) := z_i \cdot \widehat{\beta(s)}^i, \quad
\widehat{\beta(s)}^i := \frac{\beta(s)}{s - z_i}.
\]

}}

\medskip

In the second part of the paper, we illustrate the practical applicability of this theoretical framework by adapting it to examine the realization of the mammillary model describing the intravenous infusion of propofol.

\subsection{Literature review}
In the context of continuous-time linear systems, a system of the form~\eqref{sys:cont} is said to be \textit{positive} if it satisfies conditions~\eqref{eq:cond1} and~\eqref{eq:cond2} alone (see, for instance,~\cite{ de2002identification, farina2011positive, luenberger1979introduction}).
In literature, the Positive Realization Problem~\cite{benvenuti2004tutorial} is among the most studied and challenging issues in the field of positive systems.
This problem deals with identifying necessary and sufficient conditions under which a given transfer function can be realized by a positive system (see, e.g.,~\cite{farina1996minimal}, Section 3 of~\cite{farina1996existence}, Section 5 of~\cite{anderson2002nonnegative}, and~\cite{benvenuti2003minimal}).

Compartmental systems are a specific subclass within the broader category of positive systems (see, among others,~\cite{anderson1984properties, benvenuti2002model, haddad2010nonnegative, mcwilliams1985properties, walter1999compartmental}).
In particular,~\cite{maeda1977compartmental} and~\cite{benvenuti2002positive} focus on theorems that provide necessary and sufficient conditions under which a given transfer function corresponds to a compartmental system. 

Further results concerning more specific subclasses of positive systems can be found in~\cite{astolfi2004note, halmschlager2005minimal, nagy2005minimal}.
Additionally,~\cite{benvenuti2002minimal} and~\cite{benvenuti2012minimal} provide necessary and sufficient conditions for a third-order transfer function with real poles to admit a third-order (minimal) positive realization.
Further results on the positive realization problem can be found in ~\cite{benvenuti2020lower, benvenuti2020upper, benvenuti2020niep, czaja2008efficient, nagy2003lowerbound}.
For concise overviews of the established results on positive realization, see, for instance,~\cite{benvenuti2022minimal} and~\cite{rantzer2018tutorial}. 
Among various fields of application, positive realization theory has been recently applied to the study of Hessenberg forms of non-negative and Metzler matrices~\cite{grussler2022similarity}.
Moreover, it has been shown that collision avoidance within a vehicle platoon is ensured if and only if each controlled vehicle exhibits a non-negative impulse response~\cite{lunze2018adaptive},~\cite{schwab2021design}. 

This paper focuses on a particularly relevant class of mammillary systems -- a subclass of compartmental models -- in which elimination occurs exclusively from the central compartment.
These systems are characterized by the structure described in~\eqref{eq:mammillary}, where the only elimination parameter is $k_{10}$, associated solely to the central compartment.
For a comparison with mammillary systems involving elimination from all compartments, see~\cite{benet1972general}.
One may also consult~\cite{vicini2000identifiability} for a treatment of the identifiability problem in $n$-compartments linear mammillary and catenary models.
Mammillary models are frequently adapted to describe the absorption, distribution, and elimination of drugs.
For instance, they have been applied to model the pharmacokinetics of propofol, remifentanil, and rocuronium (see, e.g.,~\cite{schnider1998influence, minto1997influence, adamus2011}).
Propofol and remifentanil are intravenous anesthetic agents commonly used for the induction and maintenance of general anesthesia (see, for instance,~\cite{schnider1998influence, minto1997influence, minto1997pharmacokinetics}) accounting for hypnosis and analgesia, respectively.
Whilst rocuronium is a neuromuscular blocking agent used to facilitate endotracheal intubation and ensure muscle relaxation during surgical procedures (see, e.g.,~\cite{adamus2011}).
Additionally, mammillary models have been applied to the study of proteins metabolism kinetics in organisms, as shown in~\cite{cherruault1985general}.

\subsection{Statement of contribution}
The proposed result address the realizability of a given transfer function with a specific mammillary model. To the best of our knowledge, they are novel and have not been previously explored in literature.
The necessary and sufficient conditions for realization that we present are more restrictive compared to those in positive or compartmental realizations, due to the specific structure of the considered mammillary model.
However, despite the demand of a quite specific system structure, their formulation remains simple.

%Moreover, these results have been adapted to study the problem of realization of the mammillary model describing the infusion of propofol.
%Furthermore, the provided proofs are constructive, meaning that they explicitly detail the step-by-step procedure to obtain the realization and the parameters derivation.

\section{ Proofs of the main results}\label{proof}
Before proceeding to the proofs of Theorems \ref{thrm:conditions} and \ref{thrm:conditions_p} in the general case, we first focus on the \(3 \times 3\) case. Although simpler and more intuitive, this case remains mathematically relevant and will play an important role in the development presented in Section \ref{Appl}.

\subsection{The Three-Compartment Case}
We consider a continuous-time linear system of form~\eqref{sys:cont}, in which
\begin{equation}
\begin{gathered}
\label{eq:mammillary_3x3}
A = \left[
\begin{array}{c|ccc}
\scriptstyle -k_{10} - k_{12} - k_{13} & \scriptstyle k_{21} & \scriptstyle k_{31} \\
\hline
\scriptstyle k_{12} & \scriptstyle -k_{21} & 0 \\
\scriptstyle k_{13} & 0 & \scriptstyle -k_{31}
\end{array}
\right], \vspace{3pt}\\
B = [1, 0, 0]^T, \quad C = B^T=[1, 0, 0].
\end{gathered}
\end{equation}
Note that this system is a mammillary model of form~\eqref{eq:mammillary} consisting of three compartments.

Let $H$ be a given third order transfer function of the form $H(s)=\frac{\beta(s)}{\alpha(s)}$, where $\alpha(s)=s^3+\alpha_2 s^2+ \alpha_1 s + \alpha_0$, with $\alpha_2, \alpha_1, \alpha_0 \in \R$.
The goal is to find necessary and sufficient conditions under which $H$ has a realization of the form $H(s)=C (s I -A)^{-1} B$, where $I$ denotes the $3 \times 3$ identity matrix, and $A$, $B$, and $C$ have the specific structure of~\eqref{eq:mammillary_3x3}.

This, in turn, will allow determining coefficients $k_{10}, k_{12}, k_{13}, k_{21}, k_{23}$, so that the following identity holds: $H(s) = C(sI - A)^{-1}B$.

In our case, we have that
$C(s I -A)^{\text{adj}}B=\left(k_{21}+s\right)\,\left(k_{31}+s\right)$.

Moreover the denominator is given by
\[
\chi(s)= \det\begin{bmatrix}
k_{10}+k_{12}+k_{13}+s & -k_{21} & -k_{31}\\
-k_{12} & k_{21}+s & 0\\
-k_{13} & 0 & k_{31}+s
\end{bmatrix},
\]
that is,
\begin{align*}
\chi(s) = s^3 \!& + (k_{10} + k_{12} + k_{13} + k_{21} + k_{31})s^2 \notag \\
& + (k_{10}k_{21} + k_{13}k_{21} + k_{10}k_{31} + k_{12}k_{31} + k_{21}k_{31})s \notag \\
& + k_{10}k_{31}k_{21}
\end{align*}

 So 
\begin{equation}\label{h 3 per 3}
C (s I -A)^{-1} B=\frac{\left(k_{21}+s\right)\,\left(k_{31}+s\right)}{\chi(s)}
\end{equation}

We now state a basic algebraic lemma, accompanied by a brief proof for the sake of completeness.

\begin{lemma}
Two polynomials $p(x)$ and $q(x)$ of degree $n$ are identically equal if and only if they coincide at least in $n +1$ distinct points.
\end{lemma}

\begin{proof}
Let $r(x) =p(x)-q(x) $.
Then $r(x)$ is a polynomial of degree at most $n$.
Suppose $p(x_i) = q(x_i)$ for $n+1$ distinct points $x_1, x_2, \dots, x_{n+1}$.  
Then $r(x_i) = 0$, for all $i \in \{ 1, \dots, n+1\}$, that is, $r(x)$ has at least $n+1$ distinct zeros.

But a nonzero polynomial of degree at most $n$ can have at most $n$ distinct roots.  
Therefore, $r(x) \equiv 0$, and thus $p(x) \equiv q(x)$.

Conversely, if $p(x) \equiv q(x)$, then they clearly coincide at all points, including $n+1$ distinct ones.
\end{proof}

\begin{cor}
\label{cor:zeros}
Two monic polynomials $p(x)$ and $q(x)$ of degree $n$ are identically equal if and only if they coincide at at least $n$ distinct points.
\end{cor}

We now proceed to prove the following result in the three-dimensional case.
\begin{thrm}
\label{thm:3D}
There exist unique $k_{10}, k_{12}, k_{13}, k_{21}, k_{23}$ such that $H(s) = C (s I -A)^{-1} B$, if and only if the following conditions hold:
\begin{enumerate}
\item the relative degree of $H$ is 1;

\item the numerator $\beta$ of $H$ is a monic polynomial with simple, real, and nonzero roots.
\end{enumerate}
\end{thrm}

\begin{proof}
\begin{itemize}
\item[$(\Leftarrow)$] By hypothesis, $\beta(s)=(s-z_2)(s-z_3)$, where $z_2$ and $z_3$ are distinct, real and nonzero roots of $\beta(s)$.
Without loss of generality, we assume $z_2>z_3$.
Clearly, $H(s) = C (s I -A)^{-1} B$ if and only if $\beta(s)=(k_{21}+s)(k_{31}+s)$ and $\alpha(s) =\chi(s)$.
Equating the numerators, we get 
\begin{equation}
\label{k_21,k_31}
     k_{21}=-z_2\quad \text{ and }\quad k_{31}=-z_3. 
\end{equation}
By Corollary~\ref{cor:zeros}, since both $\alpha(s)$ and $\chi(s)$ have leading coefficient equal to 1, to ensure that $\chi(s) = \alpha(s)$, we impose the following three conditions: $\chi(0) = \alpha(0)$, $\chi(z_2) = \alpha(z_2)$, and $\chi(z_3) = \alpha(z_3)$.
From the first condition, we obtain $k_{10}=\frac{\alpha(0)}{k_{21}\,k_{31}}$,
that is, by equations~\eqref{k_21,k_31},
\begin{equation}
\label{k_10}
     k_{10}=\frac{\alpha(0)}{z_2\,z_3}.
\end{equation}
From the second condition we get $k_{12}= \frac{\alpha(z_2)}{k_{21}\,\left(k_{21}-k_{31}\right)}$, that is, again, by~\eqref{k_21,k_31},
\begin{equation}
\label{k_12}
    k_{12}= \frac{\alpha(z_2)}{z_2\,\left(z_2-z_3\right)}.
\end{equation}
From the third condition, we get $k_{13}=-\frac{\alpha(z_3)}{k_{31}\,\left(k_{21}-k_{31}\right)}$, that is, again, by~\eqref{k_21,k_31}
\begin{equation}
\label{k_13}
    k_{13}=-\frac{\alpha(z_3)}{z_3\,\left(z_2-z_3\right)}.
\end{equation}

\item[$(\Rightarrow)$] Condition~1) is trivial.
Condition~2) is necessary for equations~\eqref{k_10},~\eqref{k_12}, and~\eqref{k_13} to hold.
\end{itemize}
\end{proof}

Note that, from conditions~\eqref{k_10}--\eqref{k_13}, the five parameters of the model, $k_{10}, k_{12}, k_{13}, k_{21}$, and $k_{23}$, are uniquely determined. %up to the interchange of $k_{31} = -z_2$ and $k_{21} = -z_3$ in equations~\eqref{k_21,k_31}.

For simplicity, we present the algorithm resulting from the previous proof to compute the coefficients.

\medskip

\fbox{%
\parbox{0.95\linewidth}{%
\textit{Algorithm for computing the 5 parameters}

Let $z_2 > z_3$ be the roots of $\beta(s)$.

\medskip

\[
\begin{aligned}
k_{10} &= \frac{\alpha(0)}{z_2 z_3}, \\
k_{12} &=\frac{\alpha(z_2)}{z_2\,\left(z_2-z_3\right)}, & \quad &  k_{13}=-\frac{\alpha(z_3)}{z_3\,\left(z_2-z_3\right)},\\
k_{21} &= -z_2, & \quad & k_{31} = -z_3.\\
\end{aligned}
\]

\medskip
}}

\medskip

The following result follows directly from Theorem~\ref{thm:3D}, derived from equations~\eqref{k_21,k_31}--\eqref{k_13}.

\begin{thrm}
\label{thrm:cor3}
There exist unique $k_{10}, k_{12}, k_{13}, k_{21}, k_{23}$, all positive, such that $H(s) = C (s I -A)^{-1} B$ if and only if the following conditions hold:
\begin{enumerate}
\item the relative degree of $H$ is 1;

\item the numerator $\beta$ of $H$ is a monic polynomial with simple, real, and strictly negative roots;

\item $\alpha(0)>0$, $\alpha(z_2)<0$ and $\alpha(z_3)>0$, where $z_2$ and $z_3$ are the roots of $\beta$, with $z_2>z_3$. 
\end{enumerate}
\end{thrm}

Observe that condition~3) implies that polynomial $\alpha(s)$ has real and negative roots, in agreement with Remark~\ref{remark}.  
Indeed, from condition~3) it follows that $\alpha(s)$ admits two distinct negative real roots: one located between $0$ and $z_2$, and the other between $z_2$ and $z_3$.

Moreover, since $\alpha(s) = s^3 +\alpha_2 s^2+ \alpha_1 s + \alpha_0$, we have $\alpha(s) \to -\infty$ as $s \to -\infty$.  
Therefore, the third root of $\alpha(s)$ must also be a negative real number, lying to the left of $z_3$.

\subsection{General Case}
We shall now address the general case, namely, consider the following mammillary system of $n$-compartments
\begin{equation}
\begin{gathered}
\label{eq:mammillary_gen}
A_p = \left[
\begin{array}{c|cccc}
\scriptstyle -k_{10}-\sum_{i=2}^n k_{1i}  & \scriptstyle k_{21} & \scriptstyle k_{31} & \cdots & \scriptstyle k_{n1} \\
\hline
\scriptstyle k_{12} & \scriptstyle -k_{21} & 0 & \cdots & 0 \\
\scriptstyle k_{13} & 0 & \scriptstyle -k_{31} & \cdots & 0 \\
\vdots & \vdots & \vdots & \ddots & \vdots \\
\scriptstyle k_{1n} & 0 & 0 & \cdots & \scriptstyle -k_{n1}
\end{array}
\right], \vspace{3pt}\\
B = [1, 0, \ldots, 0]^T, \qquad C = B^T = [1, 0, \ldots, 0].
\end{gathered}
\end{equation}
where $p=(k_{10},k_{2,1},\ldots,k_{n,1},k_{1,2},\ldots,k_{1,n}) \in \Real^{2n-1}$ is the set of parameters of~\eqref{eq:mammillary_gen}.
The transfer function of~\eqref{eq:mammillary_gen} is
\[
T_p(s)=\frac{(s+k_{2,1})(s+k_{3,1})\cdots(s+k_{n,1})}{\chi_p(s)},
\]
where $\chi_p(s)=\det(sI-A_p)$.

Let $n_p(s):=(s+k_{2,1})(s+k_{3,1})\cdots(s+k_{n,1})$ be the numerator of $T_p(s)$.
    
Moreover let $n_{p,i}(s): = -k_{i1} \, \widehat{n_p(s)}^i$, where $\widehat{n_p(s)}^i$ is merely the polynomial $n_p(s)$ without the factor $(s+k_{i1})$, that is, $\widehat{n_p(s)}^i:= \frac{n_p(s)}{(s+k_{i1})}$.

With this in mind, it follows that
\begin{equation}
\label{eq:chi_p_a}
\chi_p(0)=n_p(0) k_{10}, 
\end{equation}
and
\begin{equation}
\label{eq:chi_p_b}
\chi_p(-k_{i,1})= k_{1i}\,n_{p,i}(s) \quad\text{ with } i \in \{2,\ldots,n\}.
\end{equation}
      
Now let $H(s)=\frac{\beta(s)}{\alpha(s)}$ be a given transfer function of order $n$, where $\alpha(s)=s^n+\alpha_{n-1}s^{n-1}+\alpha_{n-2}s^{n-2}+ \ldots +\alpha_1 s+ \alpha_0$, with $\alpha_i \in \R$, for $i \in \{0,1,\ldots,n-1\}$.

The generalization of Theorem~\ref{thm:3D} is the following one.
\begin{thrm}
\label{thrm:conditionss}
There exists unique $p$ such that $T_p=H$ if and only if the following conditions hold:
\begin{enumerate}
\item the relative degree of $H$ is 1;

\item the numerator of $H$ is a monic polynomial with simple, real, and nonzero roots.
\end{enumerate}
\end{thrm}
      
\begin{proof}
\begin{itemize}
\item[$(\Leftarrow)$] By assumption, the numerator of $H$ is given by
\[
\beta(s) = (s - z_2)(s - z_3) \cdots (s - z_n),
\]
that is, $\beta(s)$ has $n-1$ distinct, real, and nonzero roots.
Without loss of generality, we assume that the roots are ordered in the following way: $z_2 > z_3 > \ldots > z_n$.
Clearly, $T_p = H$ holds if and only if the numerators and denominators are equal.
So, set $k_{i1}=-z_1$, for $i \in \{2,3,\ldots,n\}$, where $z_i$ are the roots of $\beta$, that is,
\begin{equation}
\label{eq:k_i1}
k_{21}=-z_2,\ k_{31}=-z_3,\ \ldots\ k_{n1}=-z_n. 
\end{equation}

Then, the numerator of $T_p$ is equal to that of $H$.

Note that, by Corollary~\ref{cor:zeros}, $\alpha=\chi_p$ if and only if
\[
\alpha(0)=\chi_p(0), \alpha(z_i)=\chi_p(z_i), \enspace\text{for } i \in \{2,\ldots,n\}.
\]
From the first condition and by~\eqref{eq:chi_p_a}, it follows that
\begin{equation}
\label{eq:k10_1}
\beta(0) k_{10}  =\alpha(0),
\end{equation}
whose solution provides
\begin{equation}
\label{eq:k10}
k_{10}=\frac{\alpha(0)}{\beta(0)}.
\end{equation}

From the remaining conditions and by~\eqref{eq:chi_p_b}, it follows that, for $i \in \{2,\ldots,n\}$,
\begin{equation}
\label{eq:for_k1i_1}
    k_{1i} \beta_i(z_i)=\alpha(z_i),
\end{equation}
where 
\begin{equation}
\label{eq:b_i}
\beta_i(s):=z_i \widehat{\beta(s)}^i
\end{equation}
is a polynomial of degree $n-1$ , and $\widehat{\beta(s)}^i$ is $\beta(s)$ without the factor $(s-z_i)$, that is 
\[
\widehat{\beta(s)}^i:= \frac{\beta(s)}{(s-z_i)}.
\]
The solution of~\eqref{eq:for_k1i_1} provides
\begin{equation}
\label{eq:for_k1i}
     k_{1,i} =\frac{\alpha(z_i)}{\beta_i(z_i)}, \quad\text{for } i \in \{2,\ldots,n-1\}.
\end{equation}

Note that the denominators in~\eqref{eq:k10} and~\eqref{eq:for_k1i} are nonzero because of the assumptions on $\beta$.

\item[$(\Rightarrow)$]
Condition~1) is trivial.
Condition~2) is necessary to solve equations~\eqref{eq:k10_1} and~\eqref{eq:for_k1i_1}.
\end{itemize}
\end{proof}

%Note that, from conditions~\eqref{eq:k_i1},~\eqref{eq:k10} and~\eqref{eq:for_k1i}, the $2n-1$ parameters of~\eqref{eq:mammillary_gen}, that is, the components of $p$, are uniquely determined.%, up to permutations in equations~\eqref{eq:k_i1}.

The following is a direct consequence of Theorem~\ref{thrm:conditionss}, taking into account equations~\eqref{eq:k_i1},~\eqref{eq:k10} and~\eqref{eq:for_k1i}.

\begin{thrm}
There exists a unique, strictly positive $p$ such that $T_p=H$ if and only if the following conditions hold:
\begin{enumerate}
\item the relative degree of $H$ is 1;

\item the numerator of $H$ is a monic polynomial with simple, real, and strictly negative roots;

\item $H(0)>0$;

\item for any root $z_i$ of $\beta$, $\frac{\beta_i(z_i)}{\alpha(z_i)}>0$, where $\beta_i$ is defined as in~\eqref{eq:b_i}. 
\end{enumerate}
\end{thrm}

\section{Application to pharmacokinetic-pharmacodynamic models}\label{Appl}
\subsection{Model Formulation}
Anesthesia is essential in modern medicine as it ensures that patients receive surgical procedures and invasive treatments without pain, anxiety, or fear, while also protecting them from physical and psychological trauma.
Propofol is one of the most widely used hypnotic agents due to its potency, rapid redistribution, and metabolism, which ensures quick onset and short duration of action~\cite{schnider1998influence}.
Additionally, when properly dosed, it causes minimal side effects~\cite{Sahinovic2018}.

To represent intravenously administered propofol, the considered model is a pharmacokinetic-pharmacodynamic (PK-PD) one with a Wiener structure, meaning it consist of a linear component (a mammillary model) followed by a static nonlinearity (a Hill function)~\cite{schnider1998influence}. The PK describes the drug concentration trajectory inside a human body, whereas the PD studies the physiological effects of the drug, and its mechanisms of action. 
Usually, the Hill function represents the Bispectral Index (BIS), which is a numerical value derived from electroencephalogram (EEG) signals. It is commonly used in anesthesiology and it quantifies a patient’s level of consciousness, with values ranging from 100 (fully awake) to 0 (no cortical activity).

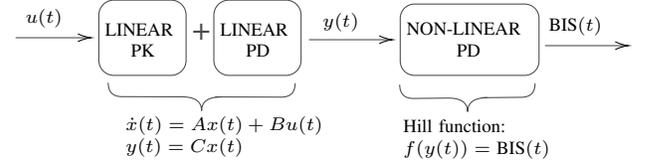
\begin{figure}[ht]
    \centering
    \tikzset{every picture/.style={line width=0.4pt}}        

    \begin{tikzpicture}[x=0.75pt,y=0.75pt,yscale=-1,xscale=1]

    % Frecce
    \draw [color={rgb, 255:red, 74; green, 74; blue, 74 }, line width=0.4pt] (157.28,95.28) -- (197.28,95.28);
    \draw [shift={(199.28,95.28)}, rotate = 180, line width=0.4pt] (6.56,-1.97) .. controls (4.17,-0.84) and (1.99,-0.18) .. (0,0) .. controls (1.99,0.18) and (4.17,0.84) .. (6.56,1.97);

    \draw [color={rgb, 255:red, 74; green, 74; blue, 74 }, line width=0.4pt] (9.28,94.28) -- (45.28,94.28);
    \draw [shift={(47.28,94.28)}, rotate = 180, line width=0.4pt] (6.56,-1.97) .. controls (4.17,-0.84) and (1.99,-0.18) .. (0,0) .. controls (1.99,0.18) and (4.17,0.84) .. (6.56,1.97);

    \draw [color={rgb, 255:red, 74; green, 74; blue, 74 }, fill=white, line width=0.4pt] (276.28,98.28) -- (316.28,98.28);
    \draw [shift={(318.28,98.28)}, rotate = 180, line width=0.4pt] (6.56,-1.97) .. controls (4.17,-0.84) and (1.99,-0.18) .. (0,0) .. controls (1.99,0.18) and (4.17,0.84) .. (6.56,1.97);

    % Parentesi graffe
    \draw [color={rgb, 255:red, 74; green, 74; blue, 74 }, line width=0.4pt] (56,117) .. controls (55.99,121.67) and (58.31,124.01) .. (62.98,124.02) -- (93.12,124.11) .. controls (99.79,124.13) and (103.11,126.47) .. (103.1,131.14) .. controls (103.11,126.47) and (106.45,124.15) .. (113.12,124.17)(110.12,124.16) -- (143.26,124.26) .. controls (147.93,124.27) and (150.27,121.95) .. (150.28,117.28);

    \draw [color={rgb, 255:red, 74; green, 74; blue, 74 }, line width=0.4pt] (203,118) .. controls (203.05,122.67) and (205.4,124.98) .. (210.07,124.93) -- (227.21,124.74) .. controls (233.88,124.67) and (237.24,126.97) .. (237.29,131.64) .. controls (237.24,126.97) and (240.54,124.6) .. (247.21,124.53)(244.21,124.56) -- (264.35,124.35) .. controls (269.02,124.3) and (271.33,121.94) .. (271.28,117.27);

    % Blocchi
    \draw [color={rgb, 255:red, 74; green, 74; blue, 74 }, line width=0.4pt] (51.28,82.88) .. controls (51.28,78.68) and (54.68,75.28) .. (58.88,75.28) -- (87.68,75.28) .. controls (91.87,75.28) and (95.28,78.68) .. (95.28,82.88) -- (95.28,105.68) .. controls (95.28,109.88) and (91.87,113.28) .. (87.68,113.28) -- (58.88,113.28) .. controls (54.68,113.28) and (51.28,109.88) .. (51.28,105.68) -- cycle;

    \draw [color={rgb, 255:red, 74; green, 74; blue, 74 }, line width=0.4pt] (109.28,82.88) .. controls (109.28,78.68) and (112.68,75.28) .. (116.88,75.28) -- (145.68,75.28) .. controls (149.87,75.28) and (153.28,78.68) .. (153.28,82.88) -- (153.28,105.68) .. controls (153.28,109.88) and (149.87,113.28) .. (145.68,113.28) -- (116.88,113.28) .. controls (112.68,113.28) and (109.28,109.88) .. (109.28,105.68) -- cycle;

    \draw [color={rgb, 255:red, 74; green, 74; blue, 74 }, line width=0.4pt] (203.28,84.48) .. controls (203.28,80.5) and (206.5,77.28) .. (210.48,77.28) -- (266.08,77.28) .. controls (270.05,77.28) and (273.28,80.5) .. (273.28,84.48) -- (273.28,106.08) .. controls (273.28,110.05) and (270.05,113.28) .. (266.08,113.28) -- (210.48,113.28) .. controls (206.5,113.28) and (203.28,110.05) .. (203.28,106.08) -- cycle;

    % Testi
    \draw (53.28,85.88) node [anchor=north west][inner sep=0.75pt] [font=\scriptsize] [align=left] {\begin{minipage}[lt]{28.51pt}\setlength\topsep{0pt}
    LINEAR
    \begin{center}
    PK
    \end{center}
    \end{minipage}};

    \draw (205.28,85.28) node [anchor=north west][inner sep=0.75pt] [font=\scriptsize] [align=left] {\begin{minipage}[lt]{46.75pt}\setlength\topsep{0pt}
    NON-LINEAR
    \begin{center}
    PD
    \end{center}
    \end{minipage}};

    \draw (97.28,85.88) node [anchor=north west][inner sep=0.75pt] [font=\scriptsize,color={rgb, 255:red, 74; green, 74; blue, 74 },opacity=1] [align=left] {{\large +}};

    \draw (13,78.4) node [anchor=north west][inner sep=0.75pt] [font=\scriptsize] {$u(t)$};

    \draw (57,132.4) node [anchor=north west][inner sep=0.75pt] [font=\scriptsize] {
    $\begin{array}{l}
    \dot{x}(t) = Ax(t) + Bu(t) \\
    y(t) = Cx(t)
    \end{array}$};

    \draw (163,80.4) node [anchor=north west][inner sep=0.75pt] [font=\scriptsize] {$y(t)$};

    \draw (197,133.4) node [anchor=north west][inner sep=0.75pt] [font=\scriptsize] {
    $\begin{array}{l}
    \text{Hill function:} \\
    f(y(t)) = \textrm{BIS}(t)
    \end{array}$};

    \draw (277,82.4) node [anchor=north west][inner sep=0.75pt] [font=\scriptsize] {$\textrm{BIS}(t)$};

    \draw (111.28,85.88) node [anchor=north west][inner sep=0.75pt] [font=\scriptsize] [align=left] {\begin{minipage}[lt]{28.51pt}\setlength\topsep{0pt}
    LINEAR
    \begin{center}
    PD
    \end{center}
    \end{minipage}};
    \end{tikzpicture}
    \caption{Schematic representation of the PK/PD Wiener model of propofol.}
    \label{fig:Wiener}
\end{figure}

For a visual representation, consider the schematic diagram in Figure~\ref{fig:Wiener}.

The aspect of the model that we focus on in this work is its linear component.
For propofol, the linear part of the Wiener model is represented by a mammillary model, more precisely, its infusion is accurately described by a three-compartment PK model of form~\eqref{eq:mammillary_3x3}, complemented by an additional effect-site compartment that captures the  PD response. 

The three compartments of the PK part are the primary, the fast, and the slow ones, each representing different physiological sites of drug distribution. The primary compartment encompasses the blood and liver, where the drug undergoes distribution and metabolism. The fast compartment corresponds to tissues with a rapid rate of drug uptake, such as muscles and viscera. The slow compartment represents tissues with a slower drug distribution, including fat and bones.
The model can be represented by the following system of differential equations:
\begin{equation}\label{sistema}
\begin{cases}
\dot{q}_1(t)\; =	&\!\!\! -(k_{10}+k_{12}+k_{13})q_1(t) \\
			&\!\!\! +k_{21}q_2(t)+k_{13}q_3(t)+u(t) \\
\dot{q}_2(t)\; =	&\!\!\! k_{12}q_1(t)-k_{21}q_2(t) \\
\dot{q}_3(t)\; =	&\!\!\! k_{13}q_1(t)-k_{31}q_3(t)\\
\dot{C}_e(t) =	&\!\!\! k_{1e}\big(\frac{q_1(t)}{V_1}\big)-k_{e0}C_e(t)
\end{cases}
\end{equation}

Input $u$ is the mass flow of the infused propofol, expressed in \unit{\milli\gram\per\second}.
Variables $q_1$, $q_2$ and $q_3$ are the drug masses, expressed in \unit{\milli\gram}, in the primary, fast and slow compartments, respectively.
Variable $C_e$ is the drug concentration in the effect-site compartment, expressed in \unit{\milli\gram\per\liter}.
The system has seven patient-dependent parameters: $k_{12}$, $k_{13}$, $k_{21}$, $k_{31}$, $k_{1e}$, $k_{10}$, and $k_{e0}$.
More precisely the drug transfer rates between compartments, expressed in \unit{\per\second}, are $k_{12}$, $k_{13}$, $k_{21}$, $k_{31}$, and $k_{1e}$.
The drug elimination rate from the primary compartment and from the effect-site compartment, expressed in \unit{\per\second}, are $k_{10}$ and $k_{e0}$, respectively.
Finally $V_1$ is the volume, expressed in \unit{\liter}, of the primary compartment.

System~\eqref{sistema} can be rewritten in the continuous-time state-space form of~\eqref{sys:cont}, in which $x(t)=\left[q_1(t), q_2(t),q_3(t), C_e(t)\right]^T$ and $y(t)=C_e(t)$, where 
\begin{equation}
\begin{gathered}
\label{A}
A =
\left[
\begin{array}{ccc|c}
\scriptstyle -k_{10}-k_{12}-k_{13} & \scriptstyle k_{21} & \scriptstyle k_{31} & \scriptstyle 0 \\
\scriptstyle k_{12} & \scriptstyle -k_{21} & \scriptstyle 0 & \scriptstyle 0 \\
\scriptstyle k_{13} & \scriptstyle 0 & \scriptstyle -k_{31} & \scriptstyle 0 \\
\hline
\scriptstyle k_{1e}/V_1 & \scriptstyle 0 & \scriptstyle 0 & \scriptstyle -k_{e0}
\end{array}
\right] \in \mathbb{R}^{4 \times 4},\vspace{3pt} \\
B = [1 , 0 , 0 , 0]^T \qquad \text{and}\qquad C = [0 , 0 , 0 , 1].
\end{gathered}
\end{equation}

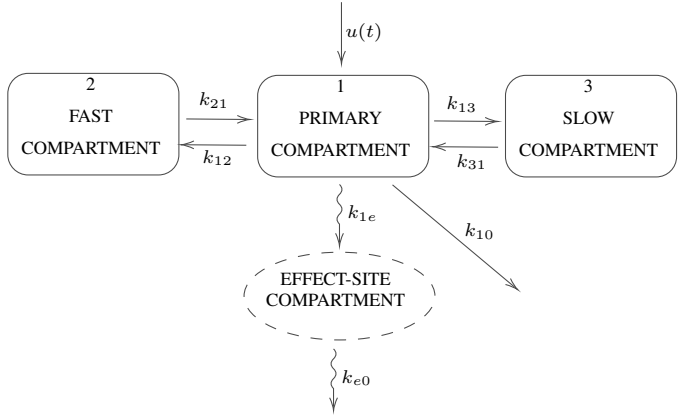
\begin{figure}[htbp]
    \centering
\tikzset{every picture/.style={line width=0.4pt}} % set default line width to 0.4pt        

\begin{tikzpicture}[x=0.75pt,y=0.75pt,yscale=-1,xscale=1]

%Rounded Rect [id:dp36244766936015493] 
\draw  [color={rgb, 255:red, 74; green, 74; blue, 74 }  ,draw opacity=1 ] (136,115.25) .. controls (136,109.59) and (140.59,105) .. (146.25,105) -- (212.02,105) .. controls (217.69,105) and (222.28,109.59) .. (222.28,115.25) -- (222.28,146.02) .. controls (222.28,151.68) and (217.69,156.27) .. (212.02,156.27) -- (146.25,156.27) .. controls (140.59,156.27) and (136,151.68) .. (136,146.02) -- cycle ;
%Rounded Rect [id:dp21102329326413938] 
\draw  [color={rgb, 255:red, 74; green, 74; blue, 74 }  ,draw opacity=1 ] (10,114.25) .. controls (10,108.59) and (14.59,104) .. (20.25,104) -- (86.02,104) .. controls (91.69,104) and (96.28,108.59) .. (96.28,114.25) -- (96.28,145.02) .. controls (96.28,150.68) and (91.69,155.27) .. (86.02,155.27) -- (20.25,155.27) .. controls (14.59,155.27) and (10,150.68) .. (10,145.02) -- cycle ;
%Rounded Rect [id:dp8285412673085883] 
\draw  [color={rgb, 255:red, 74; green, 74; blue, 74 }  ,draw opacity=1 ] (261,115.25) .. controls (261,109.59) and (265.59,105) .. (271.25,105) -- (337.02,105) .. controls (342.69,105) and (347.28,109.59) .. (347.28,115.25) -- (347.28,146.02) .. controls (347.28,151.68) and (342.69,156.27) .. (337.02,156.27) -- (271.25,156.27) .. controls (265.59,156.27) and (261,151.68) .. (261,146.02) -- cycle ;
%Straight Lines [id:da36044862886017015] 
\draw [color={rgb, 255:red, 74; green, 74; blue, 74 }  ,draw opacity=1 ]   (100,127) -- (130.28,127.26) ;
\draw [shift={(132.28,127.27)}, rotate = 180.48] [color={rgb, 255:red, 74; green, 74; blue, 74 }  ,draw opacity=1 ][line width=0.4]    (6.56,-1.97) .. controls (4.17,-0.84) and (1.99,-0.18) .. (0,0) .. controls (1.99,0.18) and (4.17,0.84) .. (6.56,1.97)   ;
%Straight Lines [id:da8251879730852576] 
\draw [color={rgb, 255:red, 74; green, 74; blue, 74 }  ,draw opacity=1 ]   (225,128) -- (255.28,128.26) ;
\draw [shift={(257.28,128.27)}, rotate = 180.48] [color={rgb, 255:red, 74; green, 74; blue, 74 }  ,draw opacity=1 ][line width=0.4]    (6.56,-1.97) .. controls (4.17,-0.84) and (1.99,-0.18) .. (0,0) .. controls (1.99,0.18) and (4.17,0.84) .. (6.56,1.97)   ;
%Straight Lines [id:da12494877224905132] 
\draw [color={rgb, 255:red, 74; green, 74; blue, 74 }  ,draw opacity=1 ]   (130.28,140.27) -- (99.28,140.03) ;
\draw [shift={(97.28,140.02)}, rotate = 0.44] [color={rgb, 255:red, 74; green, 74; blue, 74 }  ,draw opacity=1 ][line width=0.4]    (6.56,-1.97) .. controls (4.17,-0.84) and (1.99,-0.18) .. (0,0) .. controls (1.99,0.18) and (4.17,0.84) .. (6.56,1.97)   ;
%Straight Lines [id:da49777112964043524] 
\draw [color={rgb, 255:red, 74; green, 74; blue, 74 }  ,draw opacity=1 ]   (257.28,141.27) -- (226.28,141.03) ;
\draw [shift={(224.28,141.02)}, rotate = 0.44] [color={rgb, 255:red, 74; green, 74; blue, 74 }  ,draw opacity=1 ][line width=0.4]    (6.56,-1.97) .. controls (4.17,-0.84) and (1.99,-0.18) .. (0,0) .. controls (1.99,0.18) and (4.17,0.84) .. (6.56,1.97)   ;
%Straight Lines [id:da032988436947668065] 
\draw [color={rgb, 255:red, 74; green, 74; blue, 74 }  ,draw opacity=1 ]   (178.27,68.07) -- (178.28,97.27) ;
\draw [shift={(178.28,99.27)}, rotate = 269.99] [color={rgb, 255:red, 74; green, 74; blue, 74 }  ,draw opacity=1 ][line width=0.4]    (6.56,-1.97) .. controls (4.17,-0.84) and (1.99,-0.18) .. (0,0) .. controls (1.99,0.18) and (4.17,0.84) .. (6.56,1.97)   ;
%Straight Lines [id:da6499086588523354] 
\draw [color={rgb, 255:red, 74; green, 74; blue, 74 }  ,draw opacity=1 ]   (177.27,160.07) .. controls (178.94,161.74) and (178.94,163.4) .. (177.27,165.07) .. controls (175.6,166.74) and (175.6,168.4) .. (177.27,170.07) .. controls (178.94,171.74) and (178.94,173.4) .. (177.27,175.07) .. controls (175.6,176.74) and (175.6,178.4) .. (177.27,180.07) -- (177.27,181.27) -- (177.28,189.27) ;
\draw [shift={(177.28,191.27)}, rotate = 269.99] [color={rgb, 255:red, 74; green, 74; blue, 74 }  ,draw opacity=1 ][line width=0.4]    (6.56,-1.97) .. controls (4.17,-0.84) and (1.99,-0.18) .. (0,0) .. controls (1.99,0.18) and (4.17,0.84) .. (6.56,1.97)   ;
%Shape: Ellipse [id:dp2737170276300307] 
\draw  [color={rgb, 255:red, 74; green, 74; blue, 74 }  ,draw opacity=1 ][dash pattern={on 4.5pt off 4.5pt}] (129.14,216.41) .. controls (129.14,204.18) and (150.69,194.27) .. (177.28,194.27) .. controls (203.86,194.27) and (225.42,204.18) .. (225.42,216.41) .. controls (225.42,228.63) and (203.86,238.55) .. (177.28,238.55) .. controls (150.69,238.55) and (129.14,228.63) .. (129.14,216.41) -- cycle ;
%Straight Lines [id:da9480700689777702] 
\draw [color={rgb, 255:red, 74; green, 74; blue, 74 }  ,draw opacity=1 ]   (204.02,160.27) -- (265.74,211.99) ;
\draw [shift={(267.28,213.27)}, rotate = 219.96] [color={rgb, 255:red, 74; green, 74; blue, 74 }  ,draw opacity=1 ][line width=0.4]    (6.56,-1.97) .. controls (4.17,-0.84) and (1.99,-0.18) .. (0,0) .. controls (1.99,0.18) and (4.17,0.84) .. (6.56,1.97)   ;
%Straight Lines [id:da5989983952145261] 
\draw [color={rgb, 255:red, 74; green, 74; blue, 74 }  ,draw opacity=1 ]   (174.27,243.07) .. controls (175.94,244.74) and (175.94,246.4) .. (174.27,248.07) .. controls (172.6,249.74) and (172.6,251.4) .. (174.27,253.07) .. controls (175.94,254.74) and (175.94,256.4) .. (174.27,258.07) .. controls (172.6,259.74) and (172.6,261.4) .. (174.27,263.07) -- (174.27,264.27) -- (174.28,272.27) ;
\draw [shift={(174.28,274.27)}, rotate = 269.99] [color={rgb, 255:red, 74; green, 74; blue, 74 }  ,draw opacity=1 ][line width=0.4]    (6.56,-1.97) .. controls (4.17,-0.84) and (1.99,-0.18) .. (0,0) .. controls (1.99,0.18) and (4.17,0.84) .. (6.56,1.97)   ;

% Text Node
\draw (138,106) node [anchor=north west][inner sep=0.75pt]   [align=left] {\begin{minipage}[lt]{57.98pt}\setlength\topsep{0pt}
\begin{center}
{\scriptsize 1}\\{\scriptsize PRIMARY}\\{\scriptsize COMPARTMENT}
\end{center}
\end{minipage}};
% Text Node
\draw (12,105) node [anchor=north west][inner sep=0.75pt]   [align=left] {\begin{minipage}[lt]{57.98pt}\setlength\topsep{0pt}
\begin{center}
{\scriptsize 2}\\{\scriptsize FAST}\\{\scriptsize COMPARTMENT}
\end{center}
\end{minipage}};
% Text Node
\draw (263,106) node [anchor=north west][inner sep=0.75pt]   [align=left] {\begin{minipage}[lt]{57.98pt}\setlength\topsep{0pt}
\begin{center}
{\scriptsize 3}\\{\scriptsize SLOW}\\{\scriptsize COMPARTMENT}
\end{center}
\end{minipage}};
% Testo
\draw (107,141.4) node [anchor=north west][inner sep=0.75pt]  [font=\scriptsize]  {$k_{12}$};
\draw (105,112.4) node [anchor=north west][inner sep=0.75pt]  [font=\scriptsize]  {$k_{21}$};
\draw (230,114.4) node [anchor=north west][inner sep=0.75pt]  [font=\scriptsize]  {$k_{13}$};
\draw (235,143.4) node [anchor=north west][inner sep=0.75pt]  [font=\scriptsize]  {$k_{31}$};
\draw (179,77.4) node [anchor=north west][inner sep=0.75pt]  [font=\scriptsize]  {$u( t)$};
\draw (180,169.4) node [anchor=north west][inner sep=0.75pt]  [font=\scriptsize]  {$k_{1e}$};
\draw (139,203) node [anchor=north west][inner sep=0.75pt]  [font=\scriptsize] [align=center] {EFFECT-SITE\\COMPARTMENT};
\draw (239,178.4) node [anchor=north west][inner sep=0.75pt]  [font=\scriptsize]  {$k_{10}$};
\draw (177,252.4) node [anchor=north west][inner sep=0.75pt]  [font=\scriptsize]  {$k_{e0}$};

\end{tikzpicture}
    \caption{Schematic representation of the linear PK/PD compartmental model of propofol.}
    \label{fig:compartment}
\end{figure}

A schematic diagram of the linear part of the Wiener model of propofol is provided in Figure~\ref{fig:compartment}.

Note that the upper-left $3 \times 3$ submatrix of $A$ in~\eqref{A}:
\begin{equation}
\label{A_3}
A_{3\times3}:= \begin{bmatrix}
-k_a & k_{21} & k_{31}  \\
k_{12} & -k_{21} & 0 \\
k_{13} & 0 & -k_{31} 
\end{bmatrix}, 
\end{equation}
is of form~\eqref{eq:mammillary_3x3}.

A natural question that may arise is: why $k_{1e}$ does not appear among the outgoing flows from the primary compartment?
More specifically, why $a_{11}=-(k_{10}+k_{12}+k_{13})$, rather than $-(k_{10}+k_{12}+k_{13}+k_{1e})$?

The reason is that the effect-site compartment does not represent a PK compartment but rather serves as a PD extension designed to account for the lag in the effect of the drug.
Notably, the drug entering the effect-site compartment is not eliminated from the body; instead, it remains available in the blood.
Therefore, $k_ {1e}$ does not represent a loss of the drug from the primary compartment.
This is why it does not appear in $a_{11}$.

In fact, only matrix $A_{3\times 3}$ of~\eqref{A_3} characterizes the PK system consisting of three compartments: primary, fast, and slow.

\subsection{Realization of the Model}
We consider the continuous-time linear system in~\eqref {A}.
Let $H$ be a given fourth order transfer function of the form $H(s)=\frac{\beta(s)}{\alpha(s)}$, where $\alpha(s)=s^4+\alpha_3 s^3+\alpha_2 s^2+ \alpha_1 s + \alpha_0$, with $\alpha_0, \alpha_1, \alpha_2 \in \R$.

Our aim is to find necessary and sufficient conditions under which $H$, admits a realization of the form $H(s) = C(sI - A)^{-1}B$, where $I$ denotes the $4 \times 4$ identity matrix, and $A$, $B$, and $C$ have the specific structure of~\eqref {A}.
As a result, the coefficients $k_{10}, k_{12}, k_{13}, k_{21}, k_{23}, \frac{k_{1e}}{V_1}, k_{e0}$ of~\eqref{A} can be determined so that $H(s) = C(sI - A)^{-1}B$.
Observe that the following equality holds:
\[
C(sI - A)^{-1}B = \frac{\frac{k_{1e}}{V_1}(s + k_{21})(s + k_{31})}{(s + k_{e0})(s^3 + \mu_2 s^2 + \mu_1 s + \mu_0)},
\]
where
\begin{align}
\mu_2 &= k_{10} + k_{12} + k_{13} + k_{21} + k_{31}, \label{eq:mu2} \\
\mu_1 &= k_{10}k_{21} + k_{13}k_{21} + k_{10}k_{31} + k_{12}k_{31} + k_{21}k_{31}, \label{eq:mu1} \\
\mu_0 &= k_{10}k_{21}k_{31}. \label{eq:mu0}
\end{align}

%We now proceed to state and prove a theorem in this direction.

To simplify the notation, we let
\begin{gather*}
\chi(s):=(s + k_{e0})(s^3 + \mu_2 s^2 + \mu_1 s + \mu_0), \\
\Xi(s):=(s^3 + \mu_2 s^2 + \mu_1 s + \mu_0).
\end{gather*}

Notice that \( C(sI - A)^{-1}B \) is nothing but the product of equation~(\ref{h 3 per 3}) and the transfer function of a first-order low-pass filter of the form \( \frac{k_{1e}}{V_1}\frac{1}{(s + k_{e0})} \)

\begin{thrm}
\label{thrm:3D_PK}
There exist $k_{10}, k_{12}, k_{13}, k_{21}, k_{23}, \frac{k_{1e}}{V_1}, k_{e0}$ such that $H(s)=C (s I -A)^{-1} B$ if and only if the following conditions hold:
\begin{enumerate}
\item the relative degree of $H$ is 2;

\item the numerator of $H$ has simple, real, and nonzero roots;

\item the denominator of $H$ has at least one real root.
\end{enumerate}
\end{thrm}

\begin{proof}
\begin{itemize}
\item[$(\Leftarrow)$]
By hypothesis, $\beta(s)=\textit{k}(s-z_2)(s-z_3)$, where $z_2$ and $z_3$ are distinct, real and nonzero roots of $\beta(s)$.
Without loss of generality, we assume $z_2>z_3$.

Clearly, $H(s) = C (s I -A)^{-1} B$ if and only if  $\beta(s)=\frac{k_{1e}}{V_1}(k_{21}+s)(k_{31}+s)$
and $\alpha(s) =\chi(s)$.

Equating the numerators, we get 
\begin{equation}
\label{k_1e,k_21,k_31}
\frac{k_{1e}}{V_1}=\textit{k},\quad k_{21}=-z_2,\quad \text{and}\quad k_{31}=-z_3. 
\end{equation}

Let $z_0$ be the real root of $\alpha(s)$, as stated in condition~3).
This implies that $\alpha(s)$ can be written in the form $\alpha(s) = (s - z_0)\, a(s)$, where $a(s)$ is a monic polynomial of degree 3.

Equating the numerators, we get
\begin{equation}
\label{pk_e0}
k_{e0}=-z_0. 
\end{equation}

We still have to show that $a(s)=\Xi(s)$.
By Corollary~\ref{cor:zeros}, since both $a(s)$ and $\Xi(s)$ are monic polynomials, to ensure that $a(s)=\Xi(s)$, we impose, as in Theorem~\ref{thrm:conditionss} the following three conditions: $a(0)=\Xi(0)$, $a(z_2)=\Xi(z_2)$, and $a(z_3)=\Xi(z_3))$.
From the first condition, we obtain $k_{10}=\frac{a(0)}{k_{21}\,k_{31}}$, that is, by equations~\eqref{k_1e,k_21,k_31},
\begin{equation}\label{pk_10}
     k_{10}=\frac{a(0)}{z_2\,z_3}.
\end{equation}
From the second condition we get $k_{12}= \frac{a(z_2)}{k_{21}\,\left(k_{21}-k_{31}\right)}$, that is, again by~\eqref{k_1e,k_21,k_31},
\begin{equation}\label{pk_12}
    k_{12}= \frac{a(z_2)}{z_2\,\left(z_2-z_3\right)},
\end{equation}
and, from the third condition, we get $k_{13}=-\frac{a(z_3)}{k_{31}\,\left(k_{21}-k_{31}\right)}$, that is, again by~\eqref{k_1e,k_21,k_31},
\begin{equation}\label{pk_13}
    k_{13}=-\frac{a(z_3)}{z_3\,\left(z_2-z_3\right)}.
\end{equation}

\item[$(\Rightarrow)$] Condition~1) is trivial.
Condition~2) is necessary for equations~\eqref{pk_10}--\eqref{pk_13} to hold.
Condition~3) is necessary to get equation~\eqref{pk_e0}.
\end{itemize}
\end{proof}

The following result follows directly from Theorem~\ref{thrm:3D_PK}, derived from equations~\eqref{k_1e,k_21,k_31}--\eqref{pk_13}.

\begin{thrm}\label{thrm:cor3}
There exist $k_{10}, k_{12}, k_{13}, k_{21}, k_{23}, \frac{k_{1e}}{V_1}, k_{e0}$ all positive such that $H(s) = C (s I -A)^{-1} B$ if and only if the following conditions hold:
\begin{enumerate}
\item the relative degree of $H$ is 2;

\item the numerator of $H$ has simple, real, negative, and nonzero roots, and its leading coefficient is positive;

\item the denominator $\alpha(s)$ of $H(s)$ has at least one real and negative root $z_0$, that is $\alpha(s)= (s - z_0)\, a(s)$; %where $a(s)$ is a monic polynomial of degree 3.

\item $a(0)>0$, $a(z_2)<0$ and $a(z_3)>0$, where  $z_2$ and $z_3$ are the roots of $\beta$, with $z_2>z_3$. 
\end{enumerate}
\end{thrm}

Note that conditions~3) and~4) imply that denominator $\alpha(s)$ has only negative real roots.  
Thus, from condition~\eqref{pk_e0}, it can be observed that up to four different positive realizations are possible.
This depends on which of the four roots of $\alpha(s)$ is assigned to $z_0$ and consequently to filter $k_{e0}$.
Therefore, the parameters are not uniquely determined. In case of multiple positive realization, one could, for instance, choose the pole $z_0$ that is the closest to the one indicated in Schnider's model~\cite{schnider1998influence}.

%Note that, actually, if the assumed ordering  $k_{21}<k_{31}$, which has been used consistently in this work, is reversed, tat is if we modify equations~\eqref{k_1e,k_21,k_31} by setting $k_{31} = -z_2$ and $k_{21} = -z_3$, we obtain four additional possible realizations.
%Therefore, we can obtain up to eight possible realizations in total.

%In particular, in the second case, we have
%\begin{gather*}
%\frac{k_{1e}}{V_1} = \textit{k},\quad	k_{21} = -z_3,\quad	k_{31} = -z_2,\quad	 k_{e0} = -z_0, \\
%k_{10} = \frac{a(0)}{z_2\, z_3},\quad	k_{12} = -\frac{a(z_3)}{z_3\, (z_2 - z_3)},\quad	k_{13} = \frac{a(z_2)}{z_2\, (z_2 - z_3)}.
%\end{gather*}
%Hence, compared to the first case, $k_{21}$ and $k_{31}$ as well as $k_{12}$ and $k_{13}$ are interchanged.

For the sake of completeness, we briefly report the algorithm to compute the parameters of the positive realizations.

\medskip

\fbox{%
\parbox{0.95\linewidth}{%
\textit{Algorithm for computing the 7 parameters of all positive realizations}

\medskip

Let \(\beta(s)=\textit{k}(s-z_2)(s-z_3)\), with \(z_2 > z_3\).

\medskip

\[
\begin{aligned}
\frac{k_{1e}}{V_1} = \textit{k}, \quad 
k_{21} = -z_2, \quad 
k_{31} = -z_3
\end{aligned}
\]

\medskip

For any pole $z_0$ of $H(s)$,
let \(\alpha(s) = (s - z_0)\, a(s)\) and set
\[
\begin{aligned}
k_{e0} &= -z_0, \\
k_{10} &= \frac{a(0)}{z_2 z_3}, \quad
k_{12} = \frac{a(z_2)}{z_2(z_2 - z_3)}, \quad
k_{13} = -\frac{a(z_3)}{z_3(z_2 - z_3)}
\end{aligned}
\]
If $k_{10}$, $k_{12}$, $k_{13}$ are all strictly positive, keep this solution, otherwise discard it.

\medskip
}}

\medskip

Note that $k_{10}$, $k_{12}$, $k_{13}$ are all strictly positive if and only if condition 4) of Theorem~\ref{thrm:cor3} holds.

\subsection{Numerical Example}
From a practical perspective, all parameters of matrix $A$ of~\eqref{A} are strictly positive in the context of propofol infusion.  

We consider a hypothetical female patient, aged 40 years, with a height of 163~cm and a weight of 54~kg. Using these data, we apply the Schnider model~\cite{schnider1998influence} to obtain the corresponding parameters
\[
\begin{array}{l}
\frac{k_{1e}}{V_1} = 0.0018 \\
k_{21} = 0.0011 \\
k_{31} = 0.0001 \\
k_{e0} = 0.0077 \\
k_{10} = 0.0065 \\
k_{12} = 0.0063 \\
k_{13} = 0.0033
\end{array}.
\]

and the corresponding transfer function \( H(s) \), which is given by:
\[
\textstyle
H(s) = \frac{0.001792 s^2 + 2.099 \times 10^{-6} s + 1.168 \times 10^{-10}}{s^4 + 0.02482 s^3 + 0.000143 s^2 + 8.963 \times 10^{-8} s + 3.232 \times 10^{-12}}
\]

We now aim to carry out the inverse process using Theorem~\ref{thrm:cor3}: given the transfer function \( H \), we seek to recover the patient-specific parameters. As we shall see, this procedure leads us back to the parameters of the Schnider model, but not exclusively, as the inversion is not unique.
Let us proceed by verifying when $H$ satisfies conditions~1),~2),~3), and~4) therein.

Rewriting $H$ to explicitly express the roots of the numerator and the denominator, yields the following form:
\[
\textstyle
H(s) = \frac{0.0018(s+0.0001)(s+0.0011)}{(s+0.0165)(s+0.0077)(s+6.6830 \times 10^{-4})(s+3.8401\times 10^{-5})}
\]
Condition~1) is trivial.
The numerator of $H$ has simple, real, and strictly negative roots, which, following the notation of the theorem, we denote by $z_2$ and $z_3$, with $z_2 > z_3$: specifically, $z_2 = -0.0001$ and $z_3 = -0.0011$.  
Moreover, its leading coefficient is positive: $k = 0.0018$. Therefore, condition~2) is satisfied.

The denominator of $H(s)$ has real negative roots: $-0.0165$, $-0.0077$, $-6.6830 \times 10^{-4}$, and $-3.8401 \times 10^{-5}$.  
Hence, it can be expressed in the form $\alpha(s) = (s - z_0)\, a(s)$, and depending on the choice of $z_0$, we must verify if condition~4) is satisfied.

We fix $z_0 = -0.0165$.
Then, polynomial $a(s)$ is given by  
$a(s) = (s + 0.0077)(s + 6.6830 \times 10^{-4})(s + 3.8401 \times 10^{-5})$.
We observe that  
\begin{align*}
a(0)& = 1.9632 \times 10^{-10} > 0, \\
a(z_2)& = a(-0.0001) = -9.3239 \times 10^{-11} < 0, \\
a(z_3)& = a(-0.0011) = 3.1266 \times 10^{-9} > 0.
\end{align*}
Therefore, condition~4) is satisfied.
From equations~\eqref{k_1e,k_21,k_31}--\eqref{pk_13}, we then obtain:
\begin{gather*}
\textstyle
\frac{k_{1e}}{V_1} = 0.0018,\ k_{21} = -z_2 = 0.0001,\ k_{31} = -z_3 = 0.0011, \\
\textstyle
k_{e0} = -z_0 = 0.0165, \qquad k_{10} = \frac{a(0)}{z_2\, z_3} = 0.0030, \\
\textstyle
k_{12} = \frac{a(z_2)}{z_2\, (z_2 - z_3)} = 0.0015, \quad
k_{13} = -\frac{a(z_3)}{z_3\, (z_2 - z_3)} = 0.0027.
\end{gather*}
%By imposing instead $k_{21} = -z_3$ and $k_{31} = -z_2$, we obtain the following alternative realization:
%\begin{gather*}
%\textstyle
%\frac{k_{1e}}{V_1} = 0.0018,\ k_{21} = -z_3 = 0.0011,\ k_{31} = -z_2 = 0.0001, \\
%\textstyle
%k_{e0} = -z_0 = 0.0165, \qquad k_{10} = \frac{a(0)}{z_2\, z_3} = 0.0030, \\
%\textstyle
%k_{12} = -\frac{a(z_3)}{z_3\, (z_2 - z_3)} = 0.0027, \quad
%k_{13} =\frac{a(z_2)}{z_2\, (z_2 - z_3)} = 0.0015.
%\end{gather*}

Proceeding in a similar manner, by fixing $z_0 =-0.0077 $, condition~4) is again satisfied.
And the resulting values are:
\[
\begin{array}{l}
\frac{k_{1e}}{V_1} = 0.0018 \\
k_{21} = 0.0001 \\
k_{31} = 0.0011 \\
k_{e0} = 0.0077 \\
k_{10} = 0.0065 \\
k_{12} = 0.0033 \\
k_{13} = 0.0063
\end{array}.
%\quad\text{ and }\quad
%\begin{array}{l}
%\frac{k_{1e}}{V_1} = 0.0018 \\
%k_{21} = 0.0011 \\
%k_{31} =  0.0001 \\
%k_{e0} = 0.0077 \\
%k_{10} = 0.0065 \\
%k_{12} = 0.0063 \\
%k_{13} = 0.0033
%\end{array}.
\]
Note that this solution yields exactly the parameters of the Schnider model, but with \( k_{21} \) exchanged with \( k_{31} \) and \( k_{12} \) exchanged with \( k_{13} \). This arises from the initial assumption that \( k_{21} < z_{31} \).

\begin{rmrk}
Note that, actually, if the assumed ordering  $k_{21}<k_{31}$ is reversed, that is if we modify equations~\eqref{k_1e,k_21,k_31} by setting $k_{31} = -z_2$ and $k_{21} = -z_3$, we obtain 
\begin{gather*}
\frac{k_{1e}}{V_1} = \textit{k},\quad	k_{21} = -z_3,\quad	k_{31} = -z_2,\quad	 k_{e0} = -z_0, \\
k_{10} = \frac{a(0)}{z_2\, z_3},\quad	k_{12} = -\frac{a(z_3)}{z_3\, (z_2 - z_3)},\quad	k_{13} = \frac{a(z_2)}{z_2\, (z_2 - z_3)}.
\end{gather*}
Hence, compared to the initial case, $k_{21}$ and $k_{31}$ as well as $k_{12}$ and $k_{13}$ are interchanged.
\end{rmrk}

Instead, by fixing $z_0 =-6.6830 \times 10^{-4}$ or $z_0 = -3.8401 \times 10^{-5}$, condition~4) is not satisfied; therefore, a positive realization cannot be achieved. This is not so surprising since these two values are very far from the one indicated in the Schnider's model~\cite{schnider1998influence} for $k_{e0}$.

\section{Conclusions}
We studied the problem of linear systems realization through mammillary models and derived necessary and sufficient conditions under which a given transfer function admits such a realization. Compared to standard identification techniques, which may offer limited physical interpretability, mammillary models provide a structured framework aligned with the physiological mechanisms of material distribution and elimination. This property is particularly important in clinical contexts, where model interpretability and parameter meaning are essential for safe and personalized therapeutic interventions. The application to a propofol infusion model highlights the practical relevance of the proposed framework for obtaining physiologically meaningful models.

\bibliographystyle{abbrv}
\balance
\bibliography{ref}
\end{document}